\pdfoutput=1
\RequirePackage{ifpdf}
\ifpdf 
\documentclass[pdftex]{sigma}
\else
\documentclass{sigma}
\fi

\usepackage{bm}

\def\genfd{{\bm k}}

\long\def\nodo#1{{}}
\def\gg{\mathfrak{g}}

\def\RR{\mathcal{S}}

\def\hx{\hat{x}}

\def\ad{\operatorname{ad}}

\def\MR#1{} 

\begin{document}


\renewcommand{\PaperNumber}{013}

\FirstPageHeading

\ShortArticleName{Exponential Formulas and Lie Algebra Type Star Products}

\ArticleName{Exponential Formulas and Lie Algebra Type\\ Star Products}

\Author{Stjepan MELJANAC~$^\dag$, Zoran \v{S}KODA~$^\dag$ and Dragutin SVRTAN~$^\ddag$}

\AuthorNameForHeading{S.~Meljanac, Z.~\v{S}koda and D.~Svrtan}

\Address{$^\dag$~Division for Theoretical Physics, Institute
Rudjer Bo\v{s}kovi\'c,\\
\hphantom{$^\dag$}~Bijeni\v{c}ka 54, P.O.~Box 180, HR-10002
Zagreb, Croatia}
\EmailD{\href{mailto:meljanac@irb.hr}{meljanac@irb.hr}, \href{mailto:zskoda@irb.hr}{zskoda@irb.hr}}

\Address{$^\ddag$~Department of Mathematics, Faculty of Natural
Sciences and Mathematics,\\
\hphantom{$^\ddag$}~University of Zagreb, HR-10000 Zagreb, Croatia}
\EmailD{\href{mailto:dsvrtan@math.hr}{dsvrtan@math.hr}}

\ArticleDates{Received May 26, 2011, in f\/inal form March 01, 2012; Published online March 22, 2012}

\Abstract{Given formal dif\/ferential operators $F_i$ on polynomial algebra
in several variables $x_1,\ldots,x_n$, we discuss f\/inding expressions
$K_l$ determined by the equation $\exp(\sum_i x_i F_i)(\exp(\sum_j q_j x_j)) =
\exp(\sum_l K_l x_l)$ and their applications. The expressions for~$K_l$
are related to the coproducts for deformed momenta
for the noncommutative space-times
of Lie algebra type and also appear in the
computations with a class of star products.
We f\/ind combinatorial recursions
and derive formal dif\/ferential equations for f\/inding $K_l$.
We elaborate an example for a Lie algebra $su(2)$, related to a
quantum gravity application from the literature.}

\Classification{81R60; 16S30; 16S32; 16A58}

\section{Introduction}\label{section1}

Deformation quantization
\cite{ArnalCortet,ArnalCortetMolinPinczon,BayenFlato,BayenFlato2,Kathotia,kontsDefQ}
studies the associative algebra
deformations of the algebra of (usually smooth) functions
on a Poisson manifold, with a prescription
that a~linear part is proportional to the Poisson structure. Thus
the deformed algebra is, for every specialization of the
deformation parameter, isomorphic as a vector space to the
undeformed. Some deformations however can appear with a dif\/ferent
motivation, namely as algebras of functions on a space or
spacetime~\cite{AC2,BorPachol2,DimGauge,FreidelLivinePRL2006,Szabo}
underlying a noncommutative f\/ield theory (motivated by
Planck-scale physics); in that case the Poisson structure
viewpoint to deformation is not of primary concern, but rather
various additional structures and symmetries enabling to establish
elements of geometry and f\/ield theory on such a spacetime, e.g.\
the kappa-spacetime \cite{AC2,covKappaDef,scopr,MS}. Both
classical and quantum f\/ield theories can be studied on such
noncommutative spacetimes. The deformation relation to the
commutative spacetime, for such applications, is important for the
correct limit of physics at the ordinary scales, and as a {\it
method} to def\/ine various procedures by help of realizations by
commutative coordinates (and some additional operators on
commutative functions). In the deformation quantization the question of
existence and uniqueness of a~star product (possibly with additional
symmetry constraints), for a given Poisson
algebra is one of the central questions. On the other hand,
in the study of noncommutative space-time,
the noncommutative base algebra is a given at start
and exploring a distinguished isomorphism to the underlying vector space
some commutative algebra is just a method. Thus our perspective and questions
are rather dif\/ferent from the deformation quantization, though some formulas
and concepts make sense in either context.

 Technical tools to introduce geometrical notions and
calculus needed to do f\/ield theory on noncommutative spaces is in
progress. Even for simple cases, e.g.\ the case of linear Poisson
structures, the choice of a star product within an equivalence
class, can make def\/initions of additional structures more or less
accessible.

There is often the case that the generators of a noncommutative
algebra can be realized in terms of dif\/ferential operators
(elements in a Weyl algebra); sometimes one allows also formal
expressions in dif\/ferential operators up to inf\/inite order. Such
realizations are very useful in physics computations. Structures
of our concern include coproducts, deformed derivatives and
exterior calculi~\cite{allTheWay,scopr,MS,exterior} and methods
include realizations of noncommutative algebras via series in
formal dif\/ferential operators. We have been using systematically
such realizations in recent works (see e.g.~\cite{ldWeyl,covKappaDef,scopr,MS}). These realizations are also
used to treat the star products, as often an action by
noncommutative variables in dif\/ferential operators on the Fock
space gives the required isomorphism of vector spaces with a space
of commutative expressions. In general, the star product does not
extend from polynomials to the formal power series, but we can
include some subspace of formal power series; almost all
formalisms include at least the exponential series, possibly of a
polynomial argument. In particular, Fourier-type expansions of
noncommutative functions into noncommutative deformations of plane
waves, are often used~\cite{AC2,FreidelLivinePRL2006,FreidelMajid,Szabo} and many
formulas are proved in practice for bases formed by such formal
exponential expressions. Thus we here study some aspects of an
abstract version of the typical realizations of such exponentials
of formal dif\/ferential operators. We exhibit some facts and
algorithms concerning the exponential series of dif\/ferential
operators of a type often needed in this line of work and
especially in the case when we deal with a noncommutative
space-time of Lie algebra type~\cite{AC2,BorPachol2,DimGauge,MS}. Some of the mathematical
statements here may be used more generally than for the purposes
of noncommutative geometry.

 Given a f\/inite-dimensional Lie algebra $\gg$ over a f\/ield
$\genfd$, one sometimes f\/inds convenient to express the
noncommutative product on the enveloping algebra~$U(\gg)$ by
transfering it via a~vector space isomorphism to a~``Lie type star
product'' on the underlying space of the symmetric algebra~$S(\gg)$. As it is well known, both $U(\gg)$ and $S(\gg)$ have a
unique Hopf algebra structure in which the elements $x$ in $\mathfrak{g}$ are primitive, i.e.\ $\Delta(x) = 1\otimes x + x\otimes 1$. For some
purposes, e.g.\ for introducing the deformed derivatives
\cite{scopr}, the {\it coalgebra} isomorphisms $\xi :
S(\gg)\stackrel\simeq\longrightarrow U(\gg)$ (linear isomorphisms
satisfying $\Delta(\xi(x)) = (\xi\otimes\xi)(\Delta(x))$ for all
$x\in S(\gg)$), are better than other linear isomorphisms; we also
naturally require that $\xi$ restrict to the identity on $\mathbf{k}\oplus\mathfrak{g}\subset S(\mathfrak{g})$. The star product on $S(\mathfrak{g})$ transported by coalgebra isomorphism $\xi$ and def\/ined by $f*g =
\xi^{-1}(\xi(f)\cdot\xi(g))$, will then have a number of special
features, including a well-def\/ined coproduct on deformed
derivatives (interpreted as the deformed momenta~\cite{FreidelLivinePRL2006,covKappaDef,scopr,MS}) and interesting
noncommutative dif\/ferential calculi~\cite{exterior}.

The coalgebra isomorphism $\xi$ may be replaced by equivalent
data. Namely, let $\genfd$ be of characteristic zero,
in a basis $\hx_1,\ldots,\hx_n$ of $\gg$ with structure constants
$C^k_{ij}\in\genfd$ def\/ined by $[\hx_i,\hx_j]=C^k_{ij}\hx_k$ (with Einstein
summation convention); we view this basis also
as a set of generators of the enveloping algebra
$U(\gg)$. Then, by a result of $\cite{scopr}$, there is 1-1 correspondence
between coalgebra isomorphisms $\xi$'s
and matrices $(\phi^i_j)_{i,j = 1,\ldots,n}$ of formal
power series $\phi^i_j = \phi^i_j(\partial^1,\ldots,\partial^n)$
in $n$ dual variables $\partial^1,\ldots,\partial^n\in\gg^*$, such that
$\phi^i_j$'s satisfy a system of formal dif\/ferential
equations~\cite{scopr}
\[
\phi^l_j \frac{\partial}{\partial(\partial^l)}\big(\phi^k_i\big) -
\phi^l_i \frac{\partial}{\partial(\partial^l)}\big(\phi^k_j\big) =
C^s_{ij}\phi^k_s.
\]

That system of equations for $\phi^j_i$
is, on the other hand, equivalent to the requirement that the
``realization'' $\hx_i\mapsto \hx_i^\phi = \sum_{j=1}^n x_j
\phi^j_i$ extends to a homomorphism $(-)^\phi : u\mapsto u^\phi$
of associative algebras from
$U(\gg)$ into the semicompleted $n$-th Weyl algebra $A_{n,\mathbf{k}}$ (which is as a vector space identified to
$S(\mathfrak{g}\otimes\hat{S}(\mathfrak{g}^*)$; moreover as discussed in \cite{scopr}, as an
algebra, it is a smash product $S(\gg)\sharp\hat{S}(\gg^*)$).
Here $\hat{S}(V)$ denotes the completed symmetric algebra
of the f\/inite-dimensional vector space $V$
(isomorphic to the ring of formal power series in indeterminates which
form a~basis in~$V$).

While $\hx_1,\ldots,\hx_n$ do not commute, the corresponding
commutative variables $x_1,\ldots,x_n$ will be the generators of $S(\mathfrak{g})$.

 The map $\hat{x}\mapsto \hat{x}^\phi$ extends
multiplicatively to a homomorphism from $U(\gg)$ to the Weyl algebra~$\hat{A}_{n,\mathbf{k}}$ semicompleted with respect to the order of
dif\/ferential operator (we allow series in partial derivatives but
only polynomials in coordinates $x_1,\ldots, x_n$).

It has been shown in \cite{scopr} that the inverse of $\xi$ is
given by the composition of this map extending $\hat{x}\mapsto
\hat{x}^\phi$ and action of the (semicompleted) Weyl algebra on
the Fock vacuum $|0\rangle = 1_{S(\gg)}$.

 Yet another datum equivalent to $\xi$ is a vector valued
function $K = (K_1,\ldots,K_n)$ determined by the statement
\begin{gather}\label{eq:defK}
\exp\left(\sum_i k_i \hx_i^\phi\right)  \exp\left(\sum_j q_j x_j\right)  = \exp\left(\sum_l
K_l(k,q) x_l\right),
\end{gather}
where, on the left-hand side,
we use the usual Fock action of $\hat{A}_{n,\genfd}$ on
$S(\gg)$, extended {\it appropriately} to the power series
involved, and $k=(k_1,\ldots,k_n)$, $q=(q_1,\ldots,q_n)\in \mathbf{k}^n$.
If $\genfd =\mathbf{C}$ one may prefer to put $\sqrt{-1}$ in front of
all exponentials in (\ref{eq:defK}) and this introduces
Fourier-like expressions (cf.\ Section~\ref{section4}).
This article
discusses two kinds of issues:
\begin{enumerate}\itemsep=0pt
\item[--] general questions on formal operator formulas like
(\ref{eq:defK}) in formal setup.

\item[--] specif\/ics of exponentials appearing in our setup, and f\/inding
function $K$. In particular, in Section~\ref{section4}, we discuss function $K$
for a couple of realizations for $\phi$ in the case of~$su(2)$.
\end{enumerate}

The simplest case is, of course, $K(k,q) = k+q$ in which case
$\xi = e$ is the coexponential (or symmetrization) map, given by
\[
e : \ y_1\cdots y_l \mapsto \frac{1}{l!}\sum_{\sigma\in\Sigma(l)}
\hat{y}_{\sigma(1)}\cdots \hat{y}_{\sigma(l)}: \quad S(\mathfrak{g})\to U(\mathfrak{g}),
\]
for all $l$, and any elements $y_1,\ldots, y_l$ in $\mathfrak{g}\subset
S(\mathfrak{g})$, where $\hat{y}_1,\ldots,\hat{y}_l$ are the same elements as
$y_1,\ldots,y_l$, but understood as the ``noncommuting'' elements in
$\mathfrak{g}\in U(\mathfrak{g})$.

  We should f\/irst point out that in general (exception:
trivial case with $\gg$ Abelian) the star product given by $f*g =
\xi^{-1}(\xi(f)\cdot\xi(g))$ can {\it not} be continuously and
bilinearly extended to all of $\hat{S}(\gg)$; namely one can indeed f\/ind
bad power series $f$ for which $\hx_i * f = \sum_j x_j \phi^j_i(f)$ can
not be consistently written even as a formal power series (the
coef\/f\/icients of monomials in $\phi^j_i(f)$ diverge). In
particular, for an arbitrary $\phi^j_i$ it may happen that the
star product even among the exponentials of the form $\exp\big(\sum_j
q_j x_j\big)$ is ill def\/ined. To have some control of the issue, we
establish some techniques of calculating star products of
such exponentials, where one of the exponentials is
$\exp\big(\lambda \sum_j k_j x_j\big)$ and $\lambda$ is a formal commuting variable. Then the question of convergence can be studied for the result of
the calculation, after specializing our results
from formal $\lambda$ to an actual value.

While we discuss existence of
function $K = K(k,q)$ below in this article, let us suppose its
existence now and discuss the coproduct; namely if we know $K$ we
can determine the formula for coproduct (for examples see Section~\ref{section4}).

We can then introduce vector function $K^{-1}_0$ which is the
inverse of
\begin{gather}\label{eq:K0}
K_0 : \ k \mapsto K(k,0).
\end{gather} The formula for the star product for
exponentials must be a linear extension (whenever it converges) of the formula for the polynomials given by $\xi(f * g) =
\xi(f)\cdot\xi(g)$. Here $\xi$ is the inverse of the map
$U(\gg)\cong S(\gg)$ given by $u\mapsto u^\phi |0\rangle$, i.e.\ on
the monomials $\hat{x}_{j1}\cdots
\hat{x}_{jn}\mapsto\hat{x}_{j1}^\phi\cdots \hat{x}_{jn}^\phi
|0\rangle$. Thus (with $i := \sqrt{-1}$) we get
\begin{gather*}
\exp(ik\cdot x)*\exp(iq\cdot x) = \exp\big(iK^{-1}_0(k)\cdot
\hat{x}^\phi\big) \exp\big(iK^{-1}_0(q)\cdot \hat{x}^\phi\big)|0\rangle\\
\hphantom{\exp(ik\cdot x)*\exp(iq\cdot x)}{} =
\exp\big(iK^{-1}_0(k)\cdot \hat{x}^\phi\big) (\exp(iq\cdot x)).
\end{gather*}

If we def\/ine $\mathcal{D}(k,q)$ by
\[
\exp(ik\cdot x) * \exp(iq\cdot x) =
\exp(i\mathcal{D}(k,q)x),
\]
 then we obtain
\begin{gather}\label{eq:delta}
\mathcal{D}(k,q) = K\big(K^{-1}_0(k),q\big).
\end{gather}

The linear dual $U(\gg)^*$ of the enveloping algebra $U(\gg)$, is a
formal power series ring in the dual variable to the basis of $\gg$
\cite[Section~10]{ldWeyl},
with adic topology (of formal power series ring) and
equipped with a topological coproduct, namely the transpose
operator to the product~\cite{Heis}. If the dual is properly
identif\/ied with the completed symmetric algebra of the dual to the
Lie algebra $\hat{S}(\gg^*)$ with dual basis given by partial
derivatives, then this coproduct may be described by several
alternative means described in~\cite{scopr}. One of them is given
by the rule $m_*\Delta(P)(f\otimes g) = P(f * g)$ where $P\in
\hat{S}(\gg)$ is viewed as a dif\/ferential operator with constant
coef\/f\/icients evaluated at
unity and $m_*(f\otimes g) := f * g$. Then
 \begin{gather*}
\partial^j(e^{ikx}\ast e^{iqx})  = \partial^j\big( e^{i K_0^{-1}(k)\hat{x}^\phi} (e^{iqx})\big)
 =  i K_j(K_0^{-1}(k),q)e^{i K(K_0^{-1}(k),q)x}\\
 \hphantom{\partial^j(e^{ikx}\ast e^{iqx})}{}
 =  i m_\ast K_j\big(K_0^{-1}(k\otimes 1),1\otimes q\big)(e^{i k x}\otimes e^{i q x})\\
 \hphantom{\partial^j(e^{ikx}\ast e^{iqx})}{}
 =  i m_\ast K_j\big(K_0^{-1}(-i\partial\otimes 1),
1\otimes (-i\partial)\big)(e^{i k x}\otimes e^{i q x}).
\end{gather*}
Hence $\Delta (\partial^j) = i K_j(K_0^{-1}(-i\partial^1\otimes
1,\ldots,-i\partial^n\otimes 1),-i1\otimes
\partial^1,\ldots,-i1\otimes\partial^n).$

{\bf Heuristics.}
While the existence of $K_l$ satisfying (\ref{eq:defK}) is discussed
in Section~\ref{section2} in greater generality, in our case there is a simple
heuristics (pointed to the second author by S.~Skryabin whom we thank):
the exponentials in the appropriate completions of
$U(\gg)$ and $S(\gg)$ should be group-like elements (as the exponentials of
primitive elements always are) in the sense of Hopf algebras.
The coalgebra isomorphism $\xi$ preserves the property of being group-like.
We have just sketched the correspondence between the
star product and the expressions like the left-hand side of~(\ref{eq:defK}).
The product of group-like elements is group-like in a Hopf algebra;
all group-likes in~$U(\gg)$ are exponentials (in the nilpotent Lie algebra
case; otherwise we would need more precise argument for taking care of the
completions needed for the general case).

\section{Operating on exponentials}\label{section2}

This section is of a ring-theoretic nature, and deals with the facts that are more general than the study of the realizations of Lie algebra
type noncommutativity in the remainder of the article. Given a~dif\/ferential operator or a~formal power series~$F=F(d/dx)$, we
show the existence of $K = K(\lambda,k)$ such that
\[
\exp(\lambda xF(d/dx)) (\exp(kx)) = \exp(K(\lambda,k)x)
\]
as a special case of a general fact which we prove not only for the
arbitrary formal power series~$F(d/dx)$ including a multidimensional version, but even for an
arbitrary derivation~$D$ (replacing~$d/dx$), or a commuting family
of derivations, on a f\/ixed commutative ring, which is not necessarily
the polynomial or power series ring, and not necessarily in characteristic zero. Linearity of the argument of
the left-most exponential in~$x$ is, however, essential.

{\bf Basic case.} Let $A_1$ be the f\/irst Weyl algebra (over
the ring~$\mathbf{Q}$ of rational numbers or a~ring containing
$\mathbf{Q}$) with generators $x$, $d/dx$, and $k$, $\lambda$ formal commuting variables; by~$\hat{A}_1$ we will mean its completion by the
degree of dif\/ferential operator, hence allowing formal series in
$d/dx$. Then, we will show (see Corollary~\ref{sp:coralpha}) that for any $F
= F(d/dx)\in \mathbf{Q}[[d/dx]]$, there is a unique $\alpha =
\alpha(\lambda,d/dx)\in \mathbf{Q}[[\lambda,d/dx]]$ such that in
$\hat{A}_1[[\lambda,k]]$
\begin{gather}\label{eq:alphanak}
\exp(\lambda x F) = \sum_{s=0}^\infty \frac{x^s \alpha^s}{s!}.
\end{gather}
To justify~(\ref{eq:alphanak}), we will essentially use just the
fact that commuting with $x$ is a derivation of some subalgebra
containing $F$. Of course, $[\alpha(\lambda,d/dx),x]\neq 0$ in
general and the ordering between $x$ and $\alpha$ in summands $x^s
\alpha^s$ above is essential. The right-hand side may be viewed
symbolically as a normally ordered exponential $:\exp(x\alpha):$
(where $x$'s are always at the left, and $\alpha$'s always at the
right).

Regarding that $\lambda$ and $d/dx$ commute and $d/dx$ is not of f\/inite order,
we can make a substitution $d/dx\mapsto k$
and def\/ine $\alpha(\lambda, k)$ (of course the result does
not def\/ine an endomorphism of $\hat{A}_1[[\lambda,k]]$); it is a power series
with the same coef\/f\/icients but dif\/ferent argument. We will show that
\begin{gather*}
\exp(\lambda x F(d/dx))(\exp(kx)) = \exp((\alpha(\lambda,k)+k)x),
\end{gather*}
holds true in $\hat{A}_1[[\lambda,k]]$.

Similarly, in the case of several variables $x_1,\ldots,x_n$,
\begin{gather}\label{eq:alphanakmanyvar}
\exp\left(\lambda \sum_{i=1}^n x_i F_i(d/dx_1,\ldots,d/dx_n)\right) =
\sum_{s_1,\ldots,s_n}^\infty \frac{x_1^{s_1}\cdots
x_n^{s_n}\alpha_1^{s_1}\cdots \alpha_n^{s_n}}{s_1!\cdots s_n!}
\end{gather}
for unique functions $\alpha_i\in \hat{A}_n[[\lambda]]$ which, of course, depend on the commutators $[x_i, F_j]\in \hat{A}_n$.

To investigate these questions we will work in an abstract algebra
$\RR'$ which generalizes the Weyl algebra (possibly completed,
with additional parameters). By a $\mathbf{Q}$-algebra we mean an
associative unital algebra over rationals.

\begin{proposition}\label{proposition1}
Let $\mathcal{S}'$ be a $\mathbf{Q}$-algebra, and
$\mathcal{S}\subset \mathcal{S}'$ a commutative subalgebra. Let
$x\in\mathcal{S}'$ and $F\in S$ be such that the commutator $[,x]$ is
a~derivation of $\mathcal{S}$.
 Then~\eqref{eq:alphanak} is true in $\RR[[\lambda]]$ for $\alpha =
\sum_{l=0}^\infty \lambda^l A_{1,l}/l!$ $($where
$A_{1,l}\in\RR[[\lambda]]$ will be obtained below$)$.
\end{proposition}

As $\alpha$ will be explicitly constructed one evidently
has the {\it existence of}
$K(\lambda,k)=\alpha(\lambda,k) + k$ which we seek above.
Now, instead of commutator $[-,x]$ we consider an arbitrary
derivation~$D$ of~$\RR$ and we generalize the setup.

{\bf Notation.} Let $\RR$ be {\em any} commutative ring
(not necessarily containing the rationals), $F \in \RR$ an element
and $D : \RR\to \RR$ a derivation. Def\/ine a~double sequence
$\{A_{s,l}\}_{s,l\geq 0}$ of elements in $\RR$ as follows:
$A_{s,l} = 0$ unless $0\leq s\leq l$; $A_{0,0} = 1$, and
recursively $A_{s,l+1}= F\cdot (D(A_{s,l}) + A_{s-1,l})$ for
$l>1$.

{\bf Special values of $A_{r,s}$.} In particular, $A_{0,l}
= 0$ for $l>0$; $A_{1,l+1} = F D(A_{1,l})=(FD)^l(F)$ for $l\geq 0$
and $A_{s,s} = F^s$ for every $s\geq 0$.

 The evaluation of a f\/ixed derivation~$D$ of a~ring
$S$ can be represented as a~commutator with an element $T$ in an associated
extension of the original ring, the Ore extension. Its underlying
left $\mathcal{S}$-module is the free module of infinite rank underlying
the ring of polynomials in one indeterminate $\RR[T]$,
but the multiplication is changed to the unique choice which is
making it a~ring, extends $\RR\hookrightarrow \RR[T]$, and
satisf\/ies
\[
T\cdot r = r T + D(r), \qquad \forall \, r\in \RR.
\]
In particular, the Weyl algebras can be obtained as {\it iterated}
Ore extensions of polynomial rings.

\begin{theorem} \label{threcursion}  For any $\RR$,
$\RR'$, $F$, $D$ as above and $s\geq 2$, the double sequence $A_{n,l}$ satisfies the integral recursion
\begin{gather}\label{eq:Asum1sm1}
s A_{s,l} = \sum_{r=1}^{l+1-s}{l \choose r} A_{1,r} A_{s-1,l-r}.
\end{gather}
\end{theorem}

 The upper limit of the sum on the right-hand side can harmlessly
be extended up to $l-1$: the additionally included summands anyway
vanish. If we were in characteristics zero we could instead write
the recursion for the $\tilde{A}_{s,l} = s!A_{s,l}/l!$ which would
be a recursion of convolution type.

\begin{proof}[Proof of Theorem~\ref{threcursion}] The
proof is by induction on $l$: if $s>l$ the equation reads $0=0$,
for $s=l$ it reads $sA_{s,s}=sA_{s-1,s-1}F$, because $A_{s,s} = F^s$; we
just need to verify the step of induction from $(s,l)$ with $l\geq
s$ to $(s,l+1)$. For this we write $A_{s,l+1} =
FD(A_{s,l})+A_{1,1}A_{s-1,l}$, substitute~(\ref{eq:Asum1sm1}) for
$A_{s,l}$ and apply $D$ using the Leibniz rule in each summand to
obtain
\[
sA_{s,l+1} = sA_{1,1}A_{s-1,l} + \sum_{r=1}^{l-1}{l\choose r}
FD(A_{1,r})A_{s-1,l-r}+\sum_{r=1}^{l-1}{l\choose r}A_{1,r}FD(A_{s-1,l-r}).
\]
Now $FD(A_{1,r})=A_{1,r+1}$ and
\begin{gather}\label{eq:pdd}
FD(A_{s-1,l-r})=A_{s-1,l-r+1}-A_{1,1}A_{s-2,l-r},
\end{gather}
where the second summand on the right vanishes if $s=2$. Now we
f\/inish separately the case of $s=2$ and $s>2$.

For $s=2$ we obtain
\[
2A_{2,l+1} = 2A_{1,1}A_{1,l} + \sum_{r=1}^{l-1} {l\choose r}
A_{1,r+1}A_{1,l-r} + \sum_{r=1}^{l-1} {l\choose r}
A_{1,r}A_{1,l-r+1}.
\]
After absorbing $A_{1,1}A_{1,l}$ into f\/irst sum as the additional
$r=0$ summand and into the second sum as $r=l$ summand, and adding
the two sums we obtain the required form.

For $s>2$ there are several dif\/ferences. First of all
$sA_{1,1}A_{s-1,l}$ should be split into $A_{1,1}A_{s-1,l}$ which
is absorbed into the f\/irst sum as before, and
$(s-1)A_{1,1}A_{s-1,l}$ which exactly cancels the additional sum
coming from summands coming from additional $A_{1,1}A_{s-2,l-r}$
in~(\ref{eq:pdd}). The third dif\/ference is that $A_{s-1,l}A_{1,1}$
which was absorbed to extend the upper limit in the second sum for
$s=2$ does not need to be added for $s>2$ because the top limit of
$l-1$ is anyway beyond the limit of vanishing terms.
\end{proof}

\begin{corollary}\label{sp:corfactoriel} Let
$k\geq 2$ and $2\leq s = s_1 + \dots + s_k$ with $s_i \geq 1$.
Then
\[
\frac{s!}{s_1!\cdots s_k!} A_{s,l}=\sum_{l_1+\dots+l_k = l,
l>l_i\geq 1}\frac{l!}{l_1!\cdots l_k!} A_{s_1,l_1}\cdots
A_{s_k,l_k}.
\]
\end{corollary}

\begin{proof} We f\/irst prove it for $k =2$. In that case, for $s_1
= 1$ this is the statement of the theorem above. Suppose now we
have proven the statement for $s_1\geq p$. Then express $s_2 = 1 +
(s_2-1)$ and decompose $A_{s_2,l_2}$ into the sum of products of
the form $A_{1,l'_2} A_{s_3,l'_3}$, and resum~$A_{s_1,l_1}$ and~$A_{1,l'_2}$ coming from the f\/irst factor in the second sum. The
coef\/f\/icients can be easily compared.

For $k>2$ this is an easy induction on $k$ using the result for
$k=2$ both for the basis and for the step of induction.
\end{proof}

Suppose now $\RR$ is a $\mathbf{Q}$-algebra and $D$ is
$\mathbf{Q}$-linear derivation given by the commutator $[-,x]$
with a f\/ixed element $x\in\RR'$ where $\RR'\supset \RR$ is a
$\mathbf{Q}$-algebra containing $\RR$. Let $\lambda$ be a formal
variable. Then in $\RR'[[\lambda]]$

\begin{corollary}\label{sp:coralpha}
\[
\exp(\lambda x F) = \sum_{s=0}^\infty \frac{x^s \alpha^s}{s!},
\]
where $\alpha = \sum_{l=1}^\infty \lambda^l A_{1,l}/l!$ and,
of course, the commutator $[\alpha,x]\neq 0$ in general.
\end{corollary}

\begin{proof} If we set $(xF)^k = \sum_{s=1}^k x^s B_{s,k-s}$ then
we see that $B_{s,k-s}$ satisfy the recursion and initial
conditions for $A_{s,k-s}$ above. Indeed, $\sum_{s=1}^k x^s
B_{s,k-s} xF = \sum_{s=1}^k x^{s+1} B_{s,k-s} F + x^s (DB_{s,k-s})
F$ and we get the recursion after renaming the labels.

Thus the corollary follows from Corollary~\ref{sp:corfactoriel}.
\end{proof}

\begin{example}\label{p:exampleDl}
 \begin{gather*}
x \left(\frac{d}{dx}\right)^l\frac{x^{m+j(l-1)}}{(m+j(l-1))!} =
(m+j(l-1))(m+j(l-1)-1)\cdots  \\
\phantom{x \left(\frac{d}{dx}\right)^l\frac{x^{m+j(l-1)}}{(m+j(l-1))!}=}{}
  \times(m+j(l-1)-(l-1))\frac{x^{m+(j-1)(l-1)}}{(m+j(l-1))!}\\
\phantom{x \left(\frac{d}{dx}\right)^l\frac{x^{m+j(l-1)}}{(m+j(l-1))!}}{}
=   (m+(j-1)(l-1))\frac{x^{m+(j-1)(l-1)}}{(m+(j-1)(l-1))!}.
\end{gather*}
Therefore
\[
\frac{1}{j!}\left(x\frac{d^l}{dx^l}\right)^j\frac{x^{m+j(l-1)}}{(m+j(l-1))!}
= \frac{1}{j!} m(m+l-1)\cdots(m+(j-1)(l-1))\frac{x^m}{m!}.
\]
Now $(x\frac{d^l}{dx^l})^j x^n = 0$ if $m :=n-(l-1)j<0$. Therefore
\[
e^{x\frac{d^l}{dx^l}}e^{kx} =
\sum_{j=0}^\infty\sum_{m=0}^\infty\frac{(x(d/dx)^l)^j}{j!}
\frac{x^{m+j(l-1)}}{(m+j(l-1))!}k^{m+j(l-1)}.
\]
By the binomial formula
\[
\big(1-(l-1)k^{l-1}\big)^{\frac{-m}{l-1}}
 = \sum_{j=0}^\infty \frac{1}{j!} m(m-l+1)\cdots(m-(j-1)(l-1))
k^{j(l-1)}.
\]
Hence
\[
\left(\frac{k}{(1-(l-1)k^{l-1})^{\frac{1}{l-1}}}\right)^m =
\sum_{j=0}^\infty m(m-l+1)\cdots(m-(j-1)(l-1))\frac{k^{m+
j(l-1)}}{j!}.
\]
Therefore
\begin{gather}\label{eq:dxl}
e^{x\frac{d^l}{dx^l}}{e^{kx}}  =  \sum_{m=0}^\infty
\left(\frac{k}{(1-(l-1)k^{l-1})^{\frac{1}{l-1}}}\right)^m
\frac{x^m}{m!} =
\exp\left(\frac{kx}{(1-(l-1)k^{l-1})^{\frac{1}{l-1}}}\right)
 \end{gather}
for $l = 0,1,2,\ldots$. 
Therefore, for $\lambda = 1$, $F = (d/dx)^l$, we have
\[
\alpha(\lambda,k) + k = K(\lambda,k)
=\frac{k}{(1-(l-1)k^{l-1})^{\frac{1}{l-1}}}.
\]

After this work appeared at {\tt arXiv}, preprint~\cite{flajolet}
also appeared, where a formula equivalent to~(\ref{eq:dxl}) was derived as a
special case of a combinatorial method (see Fig.~1
and Chapter~6 in \cite{flajolet}).
\end{example}

 It is easy to generalize our results to treat also the
multivariable case \eqref{eq:alphanakmanyvar} via ansatz
\[
\alpha_i = \sum_{l=0}^\infty \frac{\lambda^l}{l!}
A_{0,\ldots,1,\ldots,0,l},
\] where $1$ is at $i$-th place. This
time we study a commutative algebra $\RR$ with $n$ commuting
derivations~$D_i$. The characteristics free recursion is this time
for the $(n+1)$-tuple sequence of elements
$A_{s_1,\ldots,s_n,l}\in\RR[[\lambda]]$:
\[
A_{s_1,\ldots,s_n,l+1} = \sum_{i=1}^n F_i \cdot
(D_i(A_{s_1,\ldots,s_n,l}) +
A_{s_1,\ldots,s_{i-1},(s_i)-1,s_{i+1},\ldots,s_i,l})
\]
with initial conditions $A_{0,\ldots,0,0} = 1$ and
$A_{s_1,\ldots,s_l,0}=0$ when at least one of the $s_i\neq 0$.
Then it follows by a straightforward generalization of the proof
in the case of one derivation that for all~$s_{ij}$ where $1\leq
i\leq n$, $1\leq j\leq k$ and $s_i = \sum_{j=1}^k s_{ij}$,
\[
\frac{s_1!\cdots s_k!}{s_{11}! s_{12}! \cdots s_{nk}!}
A_{s_1,\ldots,s_n,l} = \sum_{l_1+\dots+l_k = l\geq l_i\geq 1}
\frac{l!}{l_1!\cdots l_k!} A_{s_{11},\ldots,s_{1n},l_1}\cdots
A_{s_{k1},\ldots,s_{kn},l_k}.
\]

\section{Formal dif\/ferential equations}\label{section3}

We shall now exhibit some practical methods of
calculating $K(\lambda,q)$ determined
by
\[
\exp(\lambda x F(d/dx))(\exp(iqx))=\exp(K(\lambda,q)x)
\] for formal
parameter $\lambda$, real arguments $q$ and $x$, and formal series $F$.

In multivariate case, given $F(\partial) = F(\partial_1,\ldots,\partial_n)$ let
\[
K = K(\lambda, q) = (K_1(\lambda,q),\ldots,K(\lambda,q)) =
(K_1(\lambda,q_1,\ldots,q_n),\ldots,K_n(\lambda,q_1,\ldots,q_n))
\]
 be def\/ined by
\begin{gather}\label{eq:defk} e^{K(\lambda,q)\cdot x} := e^{\lambda x\cdot
F(\partial)} \left(e^{q\cdot x}\right).
\end{gather}
Then
\[
x\cdot \frac{\partial K}{\partial
\lambda}(\lambda,q)e^{K(\lambda,q) x}  =  \frac{\partial}{\partial
\lambda} \left(e^{\lambda x \cdot F(\partial)} e^{q\cdot
x}\right)  =  x\cdot F(\partial) e^{\lambda x\cdot F(\partial)}
e^{q\cdot x}.
\]
 The right-hand side can by def\/inition \eqref{eq:defk}
written as
\[
x\cdot F(\partial) e^{\lambda x\cdot F(\partial)} e^{q\cdot x} =
x\cdot F(K) e^{K(\lambda,q)\cdot x},
\]
but also as
 \begin{gather*}
e^{\lambda x\cdot F(\partial)} x\cdot F(\partial) e^{q\cdot x}  =
e^{\lambda x\cdot F(\partial)} x\cdot F(q) e^{q\cdot x}    =
\sum_{i=1}^n F_i(q) e^{\lambda x\cdot F(\partial)} x_i e^{qx} \\
\hphantom{e^{\lambda x\cdot F(\partial)} x\cdot F(\partial) e^{q\cdot x}}{}
 =  \sum_{i=1}^n F_i(q) \frac{\partial}{\partial
q_i}\big(e^{\lambda x \cdot F(\partial)} e^{q\cdot x}\big)
  =  \sum_{i=1}^n F_i(q) \frac{\partial}{\partial
q_i}\big(e^{K(\lambda,q)\cdot x}\big)  \\
\hphantom{e^{\lambda x\cdot F(\partial)} x\cdot F(\partial) e^{q\cdot x}}{}
 =  \sum_{i,j} F_i(q) \frac{\partial K_j}{\partial q_i}(\lambda,q)
x_j e^{K(\lambda,q)\cdot x}.
 \end{gather*}
Thus we obtain
\[
\sum_j x_j F_j(K) e^{K(\lambda,q)\cdot x} = \sum_{i,j} F_i(q)
\frac{\partial K_j}{\partial q_i}(\lambda,q) x_j
e^{K(\lambda,q)\cdot x}.
\]
After multiplying by $\exp(-K(\lambda,q)\cdot x)$ both sides we get
expressions linear in $x_j$. Therefore, equating the coef\/f\/icients
of $x_1,\ldots,x_n$, we obtain the system
\begin{gather*}
F_j(K(\lambda,q)) =
\sum_i F_i(q) \frac{\partial K_j}{\partial q_i}(\lambda,q) =
\frac{\partial K_j}{\partial \lambda}(\lambda,q),
\end{gather*}
where $j =1,\ldots,n$ and the boundary condition is $K(0,q) = q$.

  Let $n = 1$ and $F = (d/dx)^l$, $l>0$. Then the equations
become
\[
K^l = q^l \frac{\partial K}{\partial q} = \frac{\partial
K}{\partial\lambda},\qquad K = K(\lambda,q),\qquad K(0,q) = q.
\]
By integrating $K^l = \partial K/\partial \lambda$ we obtain that
$K^{-l+1} = (1-l)(\lambda + C(q))$ where $C = C(q)$ is some
function of~$q$. Thus
\[
\frac{\partial K^{1-l}}{\partial q} = (1-l) \frac{dC}{dq},
\]
where the left-hand side evaluates to $(1-l)K^{-l}\frac{\partial
K}{\partial q} = (1-l) K^{-l} K^l/q^l = (1-l)q^{-l}$. Therefore
$C(q) = q^{1-l}/(1-l) + C_0$ and it is easy to see that $C_0=0$.
Therefore $K^{1-l} = \lambda(1-l)+q^{1-l}$, hence, for $l>0$,
\[
e^{\lambda x\frac{d^l}{dx^l}}e^{qx} =
\exp\Big(\big(q^{1-l}+\lambda(1-l)\big)^{\frac{1}{1-l}}qx\Big) =
\exp\left(\frac{qx}{\left(1 - \lambda (l-1)
q^{l-1}\right)^{1/(l-1)}}\right),
\] in agreement with the direct
summation in Example~\ref{p:exampleDl} (for $\lambda =
1$).

{\bf A formal solution.} For a parameter $\mu$, and $1\leq i
\leq n$, def\/ine operator $Q_i(\mu)$ by
\[
Q_i(\mu) = e^{-\mu x\cdot F(\partial)}\partial_i e^{\mu x\cdot
F(\partial)} = \sum_{n=0}^\infty \frac{\mu^n\ad^n(-x\cdot
F(\partial))}{n!}(\partial_i).
\]
Now for any $R = R(\partial)$, notice
\[
[-x_j F_j(\partial), R] = F_j \frac{\partial}{\partial
(\partial_j)} R =: F_j \delta_j R,
\]
because $[F_j,R] = 0$. Then
\[
\ad^n(-x\cdot F(\partial)) F_i = - \sum_{j_1,\ldots,j_n}
F_{j_1}\delta_{j_1}(F_{j_2}\delta_{j_2}(\ldots(F_{j_n}
\delta_{j_n}(F_i))\ldots)).
\]

  Thus we obtain a {\em formal solution}
\[ Q_i(\mu) =
\partial_i + \frac{\exp(\mu \mathcal O)-1}{\mathcal O}
F_i(\partial),
\]
where $\mathcal O = \mathcal O(\partial) = \sum_i
F_i(\partial)\delta_i$. Clearly
\[
Q_i(\mu)e^{q\cdot x} = e^{-\mu x\cdot F}\partial_i e^{\mu x\cdot F}e^{q\cdot x} = e^{-\mu
x\cdot F} K_i(\mu,q) e^{K(\mu,q)\cdot x} = K_i(\mu,q).
\] Therefore
\[
K_i(\mu,q) = q_i + \frac{\exp(\mu\mathcal O(q))-1}{\mathcal O(q)}
F_i(q).
\]
For us the most important case will be $F_i(\partial) = \sum_j k_j
\phi_{ji}(\partial)$ where
\[
\sum_i x_i F_i(\partial) = \sum_{ij}
k_j x_i \phi_{ij}(\partial) = \sum_j k_j \hat{x}^\phi_j
\] for
$\hat{x}^\phi_j := \sum_i x_i\phi_{ij}(\partial)$.

The formal solution can alternatively be obtained
using the expressions $A_{r,s}$ in the recursion from Section~\ref{section2}.
Indeed, $K(\mu,q)= \alpha(\mu,q) + q$. For simplicity, we will
write it out in one variable. By Corollary~\ref{sp:coralpha},
in the notation used there, $\alpha =
\sum_{l = 1}^\infty \mu^l A_{1,l}/l!$, and the recursion gives the
special values $A_{1,l} = (F D)^{l-1} F$, for $l\geq 1$. Thus we
obtain
\[
K(\mu, q) = q + \sum_{l = 1}^\infty \mu^l (F D)^{l-1} F/l! = q +
\frac{\exp(\mu F D) - 1}{F D} F.
\]

\section[Examples related to $su(2)$]{Examples related to $\boldsymbol{su(2)}$}\label{section4}

We are now going to consider two dif\/ferent realizations
of~$su(2)$. We will slightly modify the problem: the variable
$\lambda$ will be replaced by three parameters forming a vector~$\vec{P}_1$ with length~$P_1$. General vector $q$ from above will
be denoted $\vec{P}_2$. Thus instead of $K(\lambda,q)$ we want to
f\/ind (for some realization $\phi = (\phi^a_b)$) the function $K =
K(\vec{P}_1,\vec{P}_2) = K_\phi(\vec{P}_1,\vec{P}_2)$ in the
exponent. Tricks with vector calculus and geometrically
well-chosen substitutions are useful in f\/inding the solutions.
The dif\/ferential equations will not be directly modif\/ied from the
previous section, but rather rederived on the spot in a way
introducing some useful auxiliary variables. Compare that the
formal solution from the previous section are obtained using
essentially the same variables (up to imaginary unit).

Below we shall use a basis
$\hat{x}_1$, $\hat{x}_2$, $\hat{x}_3$ of $su(2)$ satisfying
$
\lbrack\hx_a,\hx_b \rbrack = i\kappa\epsilon_{abc} \hx_c$,
where $\kappa$ is a~small parameter (this strange convention is an
adaptation for the applications to modeling some noncommutative
deformations of a space-time). Def\/ine the auxiliary variables
\[
\hat{P}_a(\mu) := e^{-i\mu k\cdot \hx}\hat{P}_a(0)e^{+i\mu k\cdot
\hx},
\]
where $\hat{P}_a(0)=\hat{p}_a = -i\partial_a$. Thus
\[
\frac{d\hat{P}_a}{d\mu}(\mu)=e^{-i\mu k\cdot
\hx}[-ik\cdot\hx,\hat{P}_a(0)]e^{+i\mu k\cdot \hx}.
\]

The realization of $U(su(2))$ of
Freidel and Livine~\cite{FreidelLivinePRL2006}.
This realization is also used in~\cite{FreidelMajid} in the context of
study of a noncommutative Fourier transform used to relate
a~group f\/ield theory related to a spin-foam model
motivated by $3d$ quantum gravity to a noncommutative f\/ield theory.
In our language their star product is coming from a realization
via formal dif\/ferential operators of inf\/inite order, is
(with Einstein summation convention) given by
\begin{gather*}
\hx_a^\phi = x_b \phi_{ba} = x_a \sqrt{1+\kappa^2 \partial^2} +
i\epsilon_{abc} \kappa x_b\partial_c,
\qquad
\phi_{ba} = \delta_{ba}
\sqrt{1+\kappa^2\partial^2}+i\kappa\epsilon_{abc}\partial_c.
\end{gather*}

Elements of $U(su(2))$ in this realization in the semicompleted
Weyl algebra act as formal dif\/ferential operators on its standard
module~-- the Fock space which is the symmetric algebra $S(su(2))$
with unit playing the role of {\it Fock vacuum} $1 = \exp(i0\cdot
x) =: | 0 \rangle$. We rescale all by imaginary units to def\/ine
$K$ by $\exp(ik\cdot\hx^\phi)\exp(iq\cdot x) = \exp(iK(k,q)\cdot
x)$. The action in the realization is
$\exp(ik\cdot\hx^\phi)|0\rangle = \exp(iK(k,0)\cdot \hx)$.

We will use the notation and the relation between $K$, $K_0$ and
the coproduct from~Section~\ref{section1}. 
For $\hat{p}_a = -i\partial_a$ we have
\begin{gather*}
[\hx_a,\hat{p}_b] = i\sqrt{1-\kappa^2 \hat{p}^2}\delta_{ab}
-i\kappa\epsilon_{abc} \hat{p}_c,
\qquad \hat{P}_a(0) = -i\partial_a.
\end{gather*}
Then $[\hat{x}_a,\partial_b] = \phi_{ab}$, what implies
\[
\frac{d\hat{P}_a}{d\mu} = k_a \sqrt{1-\kappa^2 \hat{P}^2} +
\epsilon_{abc}k_b \hat{P}_c.
\]
In these formulas the operations involving $\partial$ are
understood as acting on linear combinations of Fourier components
$\exp(iq\cdot\vec{x})$, which are the eigenvectors, with values of
$-i\partial_a$ equal to $q_a$. From now on we f\/ix a single Fourier
component $\exp(iq\cdot\vec{x})$ and write equations for $P$ which
is the corresponding eigenvalue of $\hat{P}$.

In solving the equations it is useful to utilize full vector
notation, hence writing $\vec{k}$, $\vec{q}$. We also make shortcuts
\[
L := \vec{k}\cdot\vec{P},\qquad P^2 := \sum_a (P_a)^2,
\qquad k^2 = |\vec{k}|^2, \qquad q^2 =|\vec{q}|^2.
\]
Then
\begin{gather*}
\frac{dL}{d\mu} = k^2\sqrt{1-\kappa^2 P^2},
\qquad
\frac{1}{2}\frac{dP^2}{d\mu} = L\sqrt{1-\kappa^2 P^2},
\\
\frac{d}{d\mu} \sqrt{1-\kappa^2 P^2} = -\kappa^2
\frac{dP^2/d\mu}{2\sqrt{1-\kappa^2 P^2}} =-\kappa^2 L.
\end{gather*}
Now we derive one more time,
\[
-\frac{1}{\kappa^2}\frac{d^2}{d\mu^2}\sqrt{1-\kappa^2 P^2} =
\frac{dL}{d\mu} = k^2\sqrt{1-\kappa^2 P^2}.
\]
We seek the solution of that dif\/ferential equation for
$\sqrt{1-\kappa^2 P^2}$ in the form
\[
\sqrt{1-\kappa^2 P^2} = c_1 \cos{\kappa |\vec{k}| \mu}
+ c_2 \sin{\kappa |\vec{k}| \mu}.
\]
Then of course $L = \frac{|\vec{k}|}{\kappa}(c_1 \cos{\kappa |\vec{k}| \mu} +
c_2 \sin{\kappa k \mu})$, and $ P(\mu=0) = q $, hence $c_1 =
\sqrt{1-\kappa^2 q^2}$. On the other hand, $L(\mu =
0)=\vec{q}\cdot\vec{k}$, thus  $L(\mu = 0) =|\vec{k}| \frac{c_2}{\kappa} =
\vec{q}\cdot\vec{k}$, hence $c_2 =
-\frac{\vec{q}\cdot\vec{k}}{|\vec{k}|} \kappa$. Thus
\[
L = \frac{|\vec{k}|}{\kappa}\sqrt{1-\kappa^2 q^2}\sin(\kappa|\vec{k}|\mu)+
\vec{q}\cdot\vec{k} \cos(\kappa|\vec{k}|\mu).
\]
We seek for solution for $\vec{P}$ in the form
$
\vec{P} = f_1 \vec{k} + f_2\vec{q} + f_3 \vec{k}\times \vec{q}$.
The equation
\[
\frac{d\vec{P}}{d\mu}+ \vec{k}\sqrt{1-\kappa^2 P^2} + \kappa \vec{k}\times \vec{P}
\]
becomes in these terms
\[
\frac{df_1}{d\mu}\vec{k} + \frac{df_2}{d\mu}\vec{q}+
\frac{df_3}{d\mu}\vec{k}\times\vec{q} = \sqrt{1-\kappa^2
P^2}\vec{k} + \kappa f_2 \vec{k}\times\vec{q} +
 \kappa f_3
\big((\vec{q}\cdot\vec{k})\vec{k} - k^2 \vec{q}\big),
\]
 what amounts to
the system
\begin{gather*}
\frac{df_1}{d\mu} = \vec{k}\sqrt{1-\kappa^2 P^2} + \kappa
\vec{q}\cdot\vec{k} f_3,\qquad
\frac{df_2}{d\mu} = -\kappa \vec{q}\cdot\vec{k} f_3,
\qquad
\frac{df_3}{d\mu} = \kappa f_2.
\end{gather*}
The latter two give
\[
\frac{df_2}{d\mu} = -\kappa^2 k^2 f_2,
\]
hence
\begin{gather*}
f_2 = d_1 \cos(\kappa|\vec{k}| \mu) + d_2 \sin(\kappa |\vec{k}|\mu),
\qquad
f_3 = \frac{d_1}{k}\sin(\kappa |\vec{k}|\mu) - \frac{d_2}{k}\cos(\kappa |\vec{k}|\mu).
\end{gather*}
The boundary conditions are $f_2(0) = 1$, $f_3(0) = 0$, hence $d_2
= 0$, $d_1 = 1$
\[
\vec{P} = f_1 \vec{k} + \cos(\kappa|\vec{k}|\mu) \vec{q} +
\frac{1}{|\vec{k}|}\sin(\kappa|\vec{k}|\mu) \vec{k}\times\vec{q}.
\]
Forming the inner product of this equation with $\vec{k}$ and
recalling the value of~$L$ we get the condition (both sides are
equal to~$L$)
\begin{gather*}
k^2 f_1 + \vec{q}\cdot \vec{k} \cos(\kappa|\vec{k}|\mu) =
\frac{k}{\kappa}\sqrt{1-\kappa^2 q^2} \sin(\kappa|\vec{k}| \mu) +
\vec{q}\cdot\vec{k}\cos(\kappa|\vec{k}| \mu),\\
\vec{P} = \frac{|\vec{k}|}{\kappa k}\sqrt{1-\kappa^2 q^2}\sin(\kappa
|\vec{k}| \mu) + \vec{q}\cos(\kappa |\vec{k}| \mu) + \frac{1}{|\vec{k}|}
\vec{k}\times\vec{q} \sin(\kappa |\vec{k}|\mu).
\end{gather*}
Of course, then $K(\vec{k},\vec{q}) = \vec{P}(\mu = 1)$ and
$\vec{\mathcal{D}}(\vec{k},\vec{q})$ is
then evaluated by~(\ref{eq:delta}) to obtain
\[
\vec{\mathcal{D}}(\vec{k},\vec{q}) = \sqrt{1-\kappa^2 k^2} \vec{k} +
\sqrt{1-\kappa^2 k^2}\vec{q} - \kappa\vec{k}\times\vec{q}.
\]

 The symmetric realization or ordering is def\/ined via
the condition
$
e^{i\sum_\alpha k_\alpha\hat{x}_\alpha^\phi}|0\rangle
= e^{i\sum_\alpha k_\alpha x_\alpha}$.
 In other words, $K_0$ from~(\ref{eq:K0}) is the identity. The composition of the realization
$\hat{x}\mapsto \hat{x}^\phi$ and the projection on the vacuum in
Fock space is then the inverse of the symmetrization map~\cite{ldWeyl}. We will now study $su(2)$ in this realization.

For $su(2)$ we shall now use the basis proportional to
$\sigma$-matrices $\hat{x}_i=\frac{1}{2}\sigma_i$; that basis
satisf\/ies
$[\hat{x}_i,\hat{x}_j]=i\epsilon_{ijk}\hat{x}_k$,
what follows from a useful identity $\sigma_i\sigma_j =
\delta_{ij} {\bm 1}+i\epsilon_{ijk}\sigma_k$. Then
\[
e^{ik\hat{x}} = e^{i\vec{k}\vec\sigma} = \sum_{n=0}^\infty
\frac{1}{n!}\left(i\vec{k}\frac{\vec\sigma}{2}\right)^n = \cos|k|
+ i(\vec{k}\vec\sigma)\sin|k|.
\]
In the symmetric ordering, the vector function
$\vec{\mathcal{D}}(\vec{k},\vec{q})$ from formula (\ref{eq:delta}) is
determined by
\[
e^{i\vec{q}\vec{x}} * e^{i\vec{k}\vec{x}} =
e^{i\vec{\mathcal{D}}(\vec{k},\vec{q})\vec{x}} =
\cos|\vec{\mathcal{D}}|
+\frac{i\vec{\mathcal{D}}\vec{x}}{|\vec{\mathcal{D}}|}\sin{|\vec{\mathcal{D}}|}.
\]
We need to multiply the expression in the left hand side and we
easily get
\begin{gather*}
\cos{|\vec{\mathcal{D}}|} =
\cos|\vec{k}|\cos{|\vec{q}|}-\frac{\vec{k}\vec{q}}{|\vec{k}|
|\vec{q}|} \sin{|\vec{k}|}\sin{|\vec{q}|},\\
\frac{\vec{\mathcal{D}}}{|\vec{\mathcal{D}}|}\sin{|\vec{\mathcal{D}}|}
= \frac{\vec{k}}{|\vec{k}|}\sin{|\vec{k}|}\cos{|\vec{q}|}+
\frac{\vec{q}}{|\vec{q}|}\cos{|\vec{k}|}\sin{|\vec{q}|}
-\frac{\vec{k}\times\vec{q}}{|\vec{k}||\vec{q}|}\sin{|\vec{k}|}\sin{|\vec{q}|}.
\end{gather*}
This corresponds to the realization
\[
\hat{x}_i = x_i +\frac{1}{2}\epsilon_{ijk}x_j p_k + \left(x_i
-\frac{\vec{x}\vec{p}}{p^2}p_i\right)\left(\frac{p}{2}\coth{\frac{p}{2}}-1\right),
\]
where $p_i \rightarrow -i\partial_i$. This can be used to obtain
$K$ as in the realizations above. The equation $\frac{dP_i}{d\mu}
= \phi_{ij} k_j$ is then for Fourier component
$\exp(i\vec{q}\cdot\vec{x})$
\[
\frac{dP_i}{d\mu} = k_i - \frac{1}{2}\epsilon_{ijk} k_j q_k
+\left(k_i -\frac{k_j
q_j}{q^2}q_i\right)\left(\frac{q}{2}\coth{\frac{q}{2}}-1\right).
\]
One may solve the equations looking again the solution in the form
$P(\mu) = P(\mu,\vec{k},\vec{q}) = g_1 \vec{k} + g_2 \vec{q} + g_3
\vec{k}\times \vec{q}$.

Setting $K(\vec{k},\vec{q}) = P(1,\vec{k},\vec{q})$ one obtains
$\mathcal{D}(\vec{k},\vec{q}) = K(K_0^{-1}(\vec{k},\vec{q}))$ as
before, with $K_0$ being the identity in the symmetric ordering,
hence $\mathcal{D} = K$. This way $(\vec{k},\vec{q}) \mapsto
\mathcal{D}(\mu \vec{k}, \vec{q})$ satisf\/ies the equation for $P =
P(\mu,\vec{k},\vec{q})$.

\section{Relation to works on star exponential}\label{section5}

In deformation quantization,
the $n$-th Weyl algebra $A_n$ is often identif\/ied
with a subalgebra of the Moyal algebra, i.e.~$C^\infty(V)[[h]]$,
where $V$ is the $2n$-dimensional f\/lat phase space with coordinates
$x_1,\ldots,x_n,p_1,\ldots,p_n$ with the standard symplectic form,
and the product involved is the
Moyal star product~$\star_h$~\cite{ArnalCortetMolinPinczon,BayenFlato,BayenFlato2}.
Thus the realization for generators
$\hat{x}^\phi$ can be considered as a function of the form
$\sum_i x_i\phi^i_j(p)$ in the Moyal algebra. One wants to
compute the action of $\exp\big(\sum_j k_j\hat{x}_j^\phi\big)$
on some~$g(x)$, which is
usually also an exponential.  In the Moyal representation, one replaces~$\hat{x}^\phi$ with a function of $x$ and $p$ as above and, in the power series
expansion for the exponential, replaces the usual product by the star product.
Dif\/ferential equations and other techniques for computing such star exponentials
are known (see e.g.~\cite{ArnalCortetMolinPinczon}). Now we want to act on~$g(x)$.
For this one can express the exponential involving $\hat{x}^\phi$ as
a star exponential for the Moyal star product, compute the star product
with $g(x)$ and then act on the Fock space to project
to a function of~$x$'s only.
For this, the Moyal interpretation of Weyl algebra may be suboptimal,
because the functions are in the symmetric Weyl ordering, while for the
ef\/fective computing of the action on the Fock space
(i.e.\ on the space of polynomials in $x_1,\ldots,x_n$) one usually
needs a polarized form with derivatives pushed to the right hand side.
The point of Section~\ref{section2} is essentially a method of polarizing exponentials
by a neat recursion. Our approach is instead to do the whole thing
in a single step, either by recursion
for coef\/f\/icients in formal power series (Section~\ref{section2}),
by a formal solution, 
or by solving a dif\/ferential equation.

  Our particular interest was in functions $\mathcal{D}(k,q)$ and
$K(k,q)$ from Section~\ref{section1}, 
 where it is shown that $\mathcal{D}(k,q)$
is related to the coproduct for deformed momenta.
This coproduct in the sense of
Hopf algebras, is for this case, in noncommutative
geometry interpreted as a deformed addition of momenta~\cite{Szabo}, which
is neither studied, nor has much signif\/icance
in the deformation quantization program. While in
deformation quantization one quantizes the phase space, in our situation~\cite{DimGauge,scopr,MS,Heis} one just deforms the coordinate space (thus
our star products will be just for functions of the form~$f(x)$
and not~$f(x,p)$) directly, while for f\/inding the tangent and cotangent bundles,
as well as for the deformations of the Poincar\'e algebra,
one uses Hopf algebraic techniques,
like deformed Leibniz rules~\cite{scopr}, Heisenberg double~\cite{Heis} etc.
Thus, we needed and computed very specif\/ic expressions, leading to
the computation of function $\mathcal{D}(k,q)$, hence amounting to
a new technique for computing the coproduct for deformed momenta
(for a~dif\/ferent approach see~\cite{DimGauge}).

Somewhat more generally, than in the rest of the article concerned
with Weyl algebras,
Section~\ref{section2} is concerned with certain formal expressions involving
a derivation on a general ring, and is partly beyond the scope of the usual
Weyl algebra, but the rest is about calculations involving Weyl algebras.
This is also beyond the case of Lie algebras, as such formal expressions do not
necessarily close a Lie algebra.

Realization of Lie algebras by
$\hat{x_j}^\phi = \sum_i x_i \phi^i_j(\partial)$ can be obtained
by interpreting Lie algebras as vector f\/ields on the group and computing them
in some coordinates around unit element (cf.~\cite[Sections~7--9]{ldWeyl}). More
generally, one can f\/ind similar expressions from other actions on smooth
manifolds. However, the actions do not need to exist beyond formal
neighborhood in general, as we do not ask the formal power series for~$\phi^i_j$
to have positive radius of convergence. Thus, in some cases, when the
convergence (and smoothness) allows, we can consider our expressions
as coming from a well known setup for quantization in the dif\/ferential
geometric setup. However, in full generality, the geometry of our paper
is (like in~\cite[Sections~7--9]{ldWeyl}) concerned with vector f\/ields on
formal neighborhood of the unit of the Lie group.

\section{Conclusion and further questions}\label{section6}

We have exhibited several approaches to the exponential
operators linear in variables and with arbitrary formal power
series dependence in the partial derivatives,
including direct summations, formal
operator solutions and solving dif\/ferential equations. We
have shown much detail for the case of two realizations of
$su(2)$. These equations are specif\/ically interesting
for physical applications
\cite{FreidelLivinePRL2006,Szabo,Kathotia,covKappaDef} in the
study of noncommutative spaces of Lie type via realizations by the
dif\/ferential operators of specif\/ic type.

While we def\/ined the functions $K(k,q)$ and
$\mathcal{D}(k,q)$ just formally in the relation to exponential
expressions (cf.~Section~\ref{section1}), 
computing them (up to some changes of variables)
ef\/fectively computes also the addition of momenta on the
noncommutative space, or equivalently, the coproduct on the space
of dual variables~\cite{Szabo,scopr}. This gives an important
physical application of the method present here.

We remained within a formal approach
(in the sense of formal power series). The analytic
uniformization methods from~\cite{lepowUnif} could also be used
for similar study.

One can choose some reasonably big subspace of
$\hat{S}(\gg)$ to which the star product extends well, making it a
topological algebra. Articles in deformation quantization studied
such questions also in analytic setups. But even in the simple cases,
e.g.\ when $\xi$ is the symmetrization map, def\/ining
a convenient subspace with well-def\/ined star product
and its topology is nontrivial.

The Ra\v{s}evski\v{\i}'s associative {\it hyper-envelope of a Lie algebra}
$\gg$ is a completion of $U(\gg)$ by
means of a countable family of norms $\hat{f}\mapsto
\|\hat{f}\|_{\epsilon}$ for all $\epsilon$ in an arbitrary f\/ixed family
of positive numbers having $0$ as an accumulation points, where
\[
 \|
\hat{f}\|_{\epsilon} = \max_{s_1,\ldots,
s_n}\epsilon^{-(s_1+s_2+\dots+s_n)} |f_{s_1\ldots s_n}|,
\]
for $s_1+\dots + s_N = s$, and where $f_{s_1,\ldots,s_n}$ is the
Taylor coef\/f\/icient in the front of $x_1^{s_1}\cdots x_n^{s_n}$ of the
commutative polynomial $f = e^{-1}(\hat{f})$, where $e$ is the symmetrization
map. Here $x_1,\ldots, x_n$ is
any f\/ixed basis of $\gg$, viewed as commutative coordinates. It is
nontrivial and proved by Ra\v{s}evski\v{i} in~\cite{Rasevskii}
that the algebra multiplication in $U(\gg)$
is continuous in this topology and hence that the completion of
the $U(\gg)$ as a countably normed vector space carries the unique
structure of a topological algebra extending the algebra operations on
$U(\gg)$. It may be tried to use the same def\/inition with $e$ replaced by
another coalgebra isomorphism $\xi : S(\gg)\to U(\gg)$.
The second author (Z.\v{S}.)
will show in a future publication that, under mild
conditions on $\xi$, verif\/iable in many known examples,
this modif\/ied def\/inition results in a completion of $U(\gg)$
isomorphic as a~topological algebra.

\subsection*{Acknowledgements}
We thank the Croatia MSES projects for partial supports: 098-0000000-2865 (S.M.\ and
Z.\v{S}.), 037-0372794-2807 (Z.\v{S}.) and 037-0000000-2779 (D.S.).

\pdfbookmark[1]{References}{ref}
 \LastPageEnding


\begin{thebibliography}{99}
\footnotesize\itemsep=0pt

\bibitem{AC2}
Amelino-Camelia G., Arzano M., Coproduct and star product in f\/ield theories on
  {L}ie-algebra noncommutative space-times, \href{http://dx.doi.org/10.1103/PhysRevD.65.084044}{\textit{Phys. Rev.~D}} \textbf{65}
  (2002), 084044, 8~pages, \href{http://arxiv.org/abs/hep-th/0105120}{hep-th/0105120}.

\bibitem{ArnalCortet}
Arnal D., Cortet J.C., $\ast$-products in the method of orbits for nilpotent
  groups, \href{http://dx.doi.org/10.1016/0393-0440(85)90010-5}{\textit{J.~Geom. Phys.}} \textbf{2} (1985), 83--116.

\bibitem{ArnalCortetMolinPinczon}
Arnal D., Cortet J.C., Molin P., Pinczon G., Covariance and geometrical
  invariance in $\ast$ quantization, \href{http://dx.doi.org/10.1063/1.525703}{\textit{J.~Math. Phys.}} \textbf{24}
  (1983), 276--283.

\bibitem{allTheWay}
Aschieri P., Lizzi F., Vitale P., Twisting all the way: from classical
  mechanics to quantum f\/ields, \href{http://dx.doi.org/10.1103/PhysRevD.77.025037}{\textit{Phys. Rev.~D}} \textbf{77} (2008),
  025037, 16~pages, \href{http://arxiv.org/abs/0708.3002}{arXiv:0708.3002}.

\bibitem{lepowUnif}
Barron K., Huang Y.Z., Lepowsky J., Factorization of formal exponentials and
  uniformization, \href{http://dx.doi.org/10.1006/jabr.2000.8285}{\textit{J.~Algebra}} \textbf{228} (2000), 551--579,
  \href{http://arxiv.org/abs/math.QA/9908151}{math.QA/9908151}.

\bibitem{BayenFlato}
Bayen F., Flato M., Fronsdal C., Lichnerowicz A., Sternheimer D., Deformation
  theory and quantization. I.~Deformations of symplectic structures,
  \href{http://dx.doi.org/10.1016/0003-4916(78)90224-5}{\textit{Ann. Physics}} \textbf{111} (1978), 61--110.

\bibitem{BayenFlato2}
Bayen F., Flato M., Fronsdal C., Lichnerowicz A., Sternheimer D., Deformation
  theory and quantization. II.~Physical applications, \href{http://dx.doi.org/10.1016/0003-4916(78)90225-7}{\textit{Ann. Physics}}
  \textbf{111} (1978), 111--151.

\bibitem{flajolet}
Blasiak P., Flajolet P., Combinatorial models of creation-annihilation,
  \textit{Ann. Physics} \textbf{65} (2011), Art.~B65c, 78~pages,
  \href{http://arxiv.org/abs/1010.0354}{arXiv:1010.0354}.

\bibitem{BorPachol2}
Borowiec A., Pacho{\l} A., $\kappa$-{M}inkowski spacetimes and {DSR} algebras:
  fresh look and old problems, \href{http://dx.doi.org/10.3842/SIGMA.2010.086}{\textit{SIGMA}} \textbf{6} (2010), 086, 31~pages,
  \href{http://arxiv.org/abs/1005.4429}{arXiv:1005.4429}.

\bibitem{DimGauge}
Dimitrijevi{\'c} M., Meyer F., M{\"o}ller L., Wess J., Gauge theories on the
  {$\kappa$}-{M}inkowski spacetime, \href{http://dx.doi.org/10.1140/epjc/s2004-01887-0}{\textit{Eur. Phys.~J.~C Part. Fields}}
  \textbf{36} (2004), 117--126, \href{http://arxiv.org/abs/hep-th/0310116}{hep-th/0310116}.

\bibitem{ldWeyl}
Durov N., Meljanac S., Samsarov A., {\v{S}}koda Z., A universal formula for
  representing {L}ie algebra generators as formal power series with
  coef\/f\/icients in the {W}eyl algebra, \href{http://dx.doi.org/10.1016/j.jalgebra.2006.08.025}{\textit{J.~Algebra}} \textbf{309} (2007),
  318--359, \href{http://arxiv.org/abs/math.RT/0604096}{math.RT/0604096}.

\bibitem{FreidelLivinePRL2006}
Freidel L., Livine E.R., 3D quantum gravity and ef\/fective noncommutative
  quantum f\/ield theory, \href{http://dx.doi.org/10.1103/PhysRevLett.96.221301}{\textit{Phys. Rev. Lett.}} \textbf{96} (2006), 221301,
  4~pages, \href{http://arxiv.org/abs/hep-th/0512113}{hep-th/0512113}.

\bibitem{FreidelMajid}
Freidel L., Majid S., Noncommutative harmonic analysis, sampling theory and the
  Duf\/lo map in 2+1 quantum gravity, \href{http://dx.doi.org/10.1088/0264-9381/25/4/045006}{\textit{Classical Quantum Gravity}}
  \textbf{25} (2008), 045006, 37~pages, \href{http://arxiv.org/abs/hep-th/0512113}{hep-th/0512113}.

\bibitem{Szabo}
Halliday S., Szabo R.J., Noncommutative f\/ield theory on homogeneous
  gravitational waves, \href{http://dx.doi.org/10.1088/0305-4470/39/18/030}{\textit{J.~Phys.~A: Math. Gen.}} \textbf{39} (2006),
  5189--5225, \href{http://arxiv.org/abs/hep-th/0602036}{hep-th/0602036}.

\bibitem{Kathotia}
Kathotia V., Kontsevich's universal formula for deformation quantization and
  the {C}ampbell--{B}aker--{H}ausdorf\/f formula, \href{http://dx.doi.org/10.1142/S0129167X0000026X}{\textit{Internat.~J. Math.}}
  \textbf{11} (2000), 523--551, \href{http://arxiv.org/abs/math.QA/9811174}{math.QA/9811174}.

\bibitem{kontsDefQ}
Kontsevich M., Deformation quantization of {P}oisson manifolds, \href{http://dx.doi.org/10.1023/B:MATH.0000027508.00421.bf}{\textit{Lett.
  Math. Phys.}} \textbf{66} (2003), 157--216, \mbox{\href{http://arxiv.org/abs/q-alg/9709040}{q-alg/9709040}}.

\bibitem{covKappaDef}
Meljanac S., Kre{\v{s}}i{\'c}-Juri{\'c} S., Stoji{\'c} M., Covariant
  realizations of kappa-deformed space, \href{http://dx.doi.org/10.1140/epjc/s10052-007-0285-8}{\textit{Eur. Phys.~J.~C Part. Fields}}
  \textbf{51} (2007), 229--240, \href{http://arxiv.org/abs/hep-th/0702215}{hep-th/0702215}.

\bibitem{scopr}
Meljanac S., \v{S}koda Z., Leibniz rules for enveloping algebras, \href{http://arxiv.org/abs/0711.0149}{arXiv:0711.0149}, the latest version
  available at \url{http://www.irb.hr/korisnici/zskoda/scopr5.pdf}.

\bibitem{MS}
Meljanac S., Stoji{\'c} M., New realizations of {L}ie algebra kappa-deformed
  {E}uclidean space, \href{http://dx.doi.org/10.1140/epjc/s2006-02584-8}{\textit{Eur. Phys.~J.~C Part. Fields}} \textbf{47} (2006),
  531--539, \href{http://arxiv.org/abs/hep-th/0605133}{hep-th/0605133}.



\bibitem{Rasevskii}
Ra{\v{s}}evski{\u\i} P.K., Associative superenvelope of a {L}ie algebra and its
  regular representation and ideals, \textit{Trudy Moskov. Mat. Ob\v s\v c.}
  \textbf{15} (1966), 3--54.

\bibitem{Heis}
\v{S}koda Z., Heisenberg double versus deformed derivatives,
  \href{http://dx.doi.org/10.1142/S0217751X11054772}{\textit{Internat.~J. Modern Phys. A}} \textbf{26} (2011), 4845--4854,
  \href{http://arxiv.org/abs/0909.3769}{arXiv:0909.3769}.

\bibitem{exterior}
\v{S}koda Z., Twisted exterior derivatives for enveloping algebras,
  \href{http://arxiv.org/abs/0806.0978}{arXiv:0806.0978}.


\end{thebibliography}
\end{document}